\documentclass[11pt,oneside]{amsart}

\usepackage{amscd}
\usepackage{amssymb}
\usepackage[latin1]{inputenc}
\usepackage[all]{xypic}
\usepackage{verbatim}

\def\Z{\mathbb{Z}}

\def\ZZ{\Z \oplus \Z}
\def\1{^{-1}}

\def\a{\alpha}

\def\vf{\varphi}

\def\w{\omega}
\def\G{\Gamma}

\newtheorem{thm}{Theorem}
\newtheorem{fact}{Fact}

\newtheorem{cor}{Corollary}
\newtheorem{lem}{Lemma}
\newtheorem{prop}{Proposition}

\title{Split extensions of group with infinite conjugacy classes}

\begin{document}
\maketitle
\begin{center}
{\sc Jean-Philippe PR\' EAUX}\footnote[1]{Centre de Recherche de l'Armée de l'air, Ecole de l'air, F-13661 Salon de
Provence air}\ \footnote[2]{Centre de Math\'ematiques et d'Informatique, Universit\'e de Provence, 39 rue
F.Joliot-Curie, F-13453 marseille
cedex 13\\
\indent {\it E-mail :} \ preaux@cmi.univ-mrs.fr}
\end{center}

\begin{abstract}
We give a characterization of the group property of being with infinite conjugacy classes (or {\it icc}, {\it i.e.}
$\not= 1$ and of which all
conjugacy classes beside 1 are  infinite) for split extensions of group.\\
\end{abstract}

\section*{Introduction}

A group is said to be with {\sl infinite conjugacy classes} (or {\it icc}) if it is non trivial, and if all its
conjugacy classes beside $\{ 1\}$ are infinite. This property is motivated by the theory of Von Neumann algebra, since
for any group $\G$, a necessary and sufficient condition for its Von Neumann algebra $W_\lambda^*(\G)$ to be a type
$II-1$ factor is that $\G$ be icc (cf. \cite{roiv}).

The property of being icc has been characterized in several classes of groups : 3-manifolds and $PD(3)$ groups in
\cite{aogf3v}, groups acting on Bass-Serre trees in \cite{ydc}, wreath products and finite extensions in \cite{p1,p2}.
We will focus here on groups defined by a split extension (also called {\sl semi-direct product}).\medskip\\
\indent
 Towards this direction particular results are already known. In \cite{p2} has been proved the following :
:\medskip\\
{\it \indent Let $G$ be a finite extension of $K$ :
$$1\longrightarrow K\longrightarrow G\longrightarrow \text{$Q$ finite}\longrightarrow 1$$
then $G$ is icc if and only if $K$ is icc and the natural homomorphism $\theta : Q\longrightarrow Out(K)$ is
injective.}
\medskip\\
In particular, it applies when the finite extension splits :\medskip\\
{\it  \indent Let $G=K\rtimes_\theta Q$, with $Q$ finite ; $G$ is icc if and only if $K$ is icc and the homomorphism
$Q\longrightarrow
Out(K)$ induced by $\theta$ is injective.}\bigskip\\
\indent In \cite{ydc},  has been proved, among other results, the following characterization of icc extensions by $\Z$
 :\smallskip\\
\indent
 {\it Let $G=D\rtimes_\theta \Z$ with $D\not=\{1\}$ ; then $G$ is
not icc if and only if one of the following conditions is satisfied :
\begin{itemize}
\item[(i)] D contains a $\theta(Q)$-stable normal subgroup $N\not=\{1\}$ and either $N$ is finite or $D=\Z^n$ and if
$\pi$ is the natural homomorphism from $G$ to $GL(n,\Z)$ extending $\theta$, then $N$ has only finite $\pi(G)$-orbits,
\item[(ii)] the homomorphism $\Z\longrightarrow Out(D)$ induced by $\theta$ is non injective.
\end{itemize}
}\bigskip

We give a generalization of these partial results by proposing a general characterization of split extensions with
infinite conjugacy classes.

\section{Preliminaries} Let $G$ be a group, $H$ a non empty subset of $G$,  and $u,g$  elements of $G$ ; then
$Z_G(H)$ and $Z(G)$ denote respectively the centralizer of $H$ in $G$ and  the center of $G$. The element $u^g$ of $G$
is defined as $u^g=g^{-1}ug$, while $u^H=\{u^g\ |\ g\in H\}$ ; in particular $u^G$ denotes the conjugacy class of $u$
in $G$. One immediatly verifies that the cardinality of $u^G$ equals the index of $Z_G(u)$ in $G$ so that $u^G$ is
finite if and only if $Z_G(u)$ has a finite index in $G$.

The set of elements having a finite conjugacy class in $G$ turns out to be a characteristic subgroup of $G$ that we
denote by $FC(G)$ ; it is a so called $FC$-group, that is a group whom all conjugacy classes are finite. Obviously $G$
is icc if and only if $FC(G)=\{1\}$. The class of $FC$-groups has been extensively studied and it's a well known fact
that finitely generated $FC$-groups are precisely those groups defined by a central finite extension of a f.g. abelian
group (c.f. \cite{fc}). In particular, in any f.g. $FC$-groups $K$, the subset of torsion elements $Tor(K)$ is a
characteristic subgroup of $K$ and the quotient $K/Tor(K)$ is free abelian with a finite rank.\smallskip\\
\indent In the following the group $G$ stands for the split extension $G=K\rtimes_\theta Q$ (or {semi-direct product})
with normal factor $K$, retract factor $Q$ and associated homomorphism $\theta: Q\longrightarrow Aut(K)$ ; with these
notations, for any $k\in K$, $q\in Q$, $q^{-1}kq=\theta(q)(k)$. Let $\pi : G\longrightarrow Aut(K)$ be  the
homomorphism defined by $\forall\, g\in G,\, k\in K$, $\pi(g)(k)=g^{-1}kg$ ; it extends on $G$ both $\theta$ and the
natural homomorphism $\pi_K : K\longrightarrow Inn(K)$, that is the diagram below commutes.
\[\xymatrix{
1\ar[r]&K\ar[r]\ar_{\pi_{K}}[d] & G \ar[r] \ar_{\pi}[d] & Q \ar[r]\ar^{\theta}[d]& 1\\
1 \ar[r] &Inn(K) \ar[r]&\pi(G) \ar[r]& \theta(Q)
\ar[r]&1\\
 }
\]
The subgroup $\pi(G)$ of $Aut(K)$ is an extension of $Inn(K)$ by $\theta(Q)$ ; in general the extension does not split,
despite the above one does.

We shall write in the following $\theta_q$ and $\pi_g$ instead of $\theta(q)$ and $\pi(g)$. We will denote by ${\Theta}
: FC(Q)\longrightarrow Out(K)$ the homomorphism induced by $\theta:Q\longrightarrow Aut(Q)$.

\section{Statement of the main result}

The first result we prove is the following characterization of semi-direct products with infinite conjugacy classes :

\begin{thm}\label{semidirect}
Let $G=K\rtimes_\theta Q\not=1$ be a split extension,  and $\pi : G\longrightarrow Aut(K)$,\linebreak[5]
$\Theta:FC(Q)\longrightarrow Out(K)$ be the homomorphisms defined as above.
%
Then $G$ is not icc if and only if one of the following conditions is satisfied :
\begin{itemize}
\item[$(i)$] $K$  contains a normal subgroup $N\not=1$ preserved under the action of $\pi(G)$ and such that either $N$
is finite, or $N\approx \Z^n$ has only finite $\pi(G)$-orbits.\smallskip
 \item[$(ii)$]
 $\ker \Theta$ contains $q\not=1$ with
 $\forall\,x\in K,\,\theta_q(x)=k^{-1}xk$,
 for some $k\in K$ with finite $\theta(Q)$-orbit.
\end{itemize}
\end{thm}

\noindent{\bf Remark 1.} Condition $(ii)$ can be rephrased as :
\begin{itemize}{\it
 \item[$(ii)$] either $\theta :FC(Q)\longrightarrow Aut(K)$ is non injective or
 the subgroup $\pi_K^{-1}(\theta(FC(Q)))$ of $K$ contains  $k\not=1$ whose $\theta(Q)$-orbit is finite.}
\end{itemize}
\noindent {\bf Example.} Suppose $G=K\rtimes_\theta Q$ ; if $K$ satisfies any of the above assumptions, then $G$ is not
icc :\smallskip\\
-- $K$ is a non trivial elementary group,\smallskip\\
-- $Z(K)$ contains a non trivial finite subgroup,\smallskip\\
-- $FC(K)\setminus 1$ contains a finite $\theta(Q)$-orbit,\smallskip\\
-- $Tor(FC(K))$ is a non trivial finite group.\smallskip\\
(in each case condition $(i)$ of theorem \ref{semidirect} is satisfied.)\\
-- $\theta$ is non injective,\\ (in which case condition $(ii)$ follows.)\smallskip

The theorem \ref{semidirect} can be rephrased in several ways. The first rephrasing is by mean of the finite
$\theta(Q)$-orbits in $K$.
\smallskip

\noindent{\bf Theorem 2.}\ {\it Let $O_\theta$ be the union of all finite $\theta(Q)$-orbits in $K$ ; $G$ is icc if and
only if :\begin{itemize}\item[$(a)$] $O_\theta\cap FC(K)=1$, and\item[$(b)$] $O_\theta\cap
\pi_{K}^{-1}(\theta(FC(Q)))=1$, and \item[$(c)$] the restricted homomorphism $\theta : FC(Q)\longrightarrow Aut(K)$ is
injective.
\end{itemize}
}\smallskip

Condition $(a)$ of theorem 2 is equivalent to the negation of condition $(i)$ of theorem \ref{semidirect}. Negation of
condition $(ii)$ of theorem \ref{semidirect} is equivalent to the conjunction of conditions $(b)$ and $(c)$ of theorem
2. So that theorem 2 can be seen as a way of reducing condition $(ii)$ into the obvious condition :\ {\it $\theta :
FC(Q)\longrightarrow Aut(K)$ non injective}, and a residual one.

In this direction one can also reduce condition $(ii)$ of theorem \ref{semidirect} into the condition that {\it either
$\theta :FC(Q)\longrightarrow Aut(K)$ is non injective or $G$ contains --what we called-- a twin
$FC$-subfactor}.\smallskip\\
\noindent{\bf Theorem 3.}\ {\it In theorem \ref{semidirect}, condition $(ii)$ can be changed into :
\begin{itemize}
\item[$(ii.a)$] $\theta:FC(Q)\longrightarrow Aut(K)$ is non injective, or \item[$(ii.b)$] $G$ contains a twin
$FC$-subfactor $C\rtimes C$.
\end{itemize}
}

 Roughly speaking, a twin $FC$-subfactor is  a transversal subgroup, either
 $\Z^n\times \Z^n$ or $C\rtimes C$ for $C$ a finite group, which is $\theta(Q)$-stable with only finite
 $\theta(Q)$-orbits, and such that
$\theta^{-1}\circ \pi$ sends isomorphically the left factor on the right one. (cf. \S 6).\\

\section{Proof of theorem \ref{semidirect}}

This section is entirely devoted to proving theorem \ref{semidirect}.\\

 \noindent{\it Proof of theorem \ref{semidirect}.}

 We first prove the sufficient part of the assumption, that is,  if either condition
$(i)$ or condition $(ii)$ is satisfied, then $G$ is not icc.
\begin{fact} Condition {$(i)$} implies that $G$ is not icc.\end{fact}
\noindent{\it Proof of the fact 1.} Suppose the condition $(i)$ is satisfied. Since the conjugacy class in $G$ of an
element of $K$ is its orbit under the action of $\pi(G)$, obviously each element of $N$  has a finite conjugacy class
in $G$, and hence $G$ is not icc.\hfill$\square$
\begin{fact} Condition {$(ii)$} implies that $G$ is not icc.\end{fact}
\noindent{\it Proof of the fact 2.} Suppose the condition $(ii)$ is satisfied ; let $\w=k^{-1}q\not=1$, so that
$Z_G(\w)\supset K$. Let $Stab_\theta(k)$ denotes the stabilizer of $k$ in $\theta(Q)$ ; since it has a finite index in
$\theta(Q)$, then $Q_0=\theta^{-1}(Stab_\theta(k))$ has a finite index in $Q$. Hence $Q_1=Z_Q(q)\cap Q_0$ also has a
finite index in $Q$. Then for any $u\in Q_1$, $\w^u=u^{-1}k^{-1}qu=\theta_u(k^{-1})q^u=k^{-1}q=\w$ ; hence
$Z_G(\w)\supset Q_1$. It follows that $Z_G(\w)$ contains $K\rtimes Q_1$ and hence has a finite index in $G$, so that
$G$ is not icc.\hfill$\square$\smallskip
We now prove the necessary part of the assumption, that is, if $G$ is not icc then either condition $(i)$ or condition
$(ii)$ is satisfied. Let $G$ be not icc : since $G\not=1$, there exists $u\not=1$ in $G$ such that $u^G$ is finite.
\begin{fact} If $K$ contains $u\not=1$ with $u^G$ finite, then condition $(i)$ follows.\end{fact}
\noindent{\it Proof of the fact 3.}
 Let $N'$ be the subgroup
of $K$ finitely generated by the set $u^G$. Then $N'$ is
 preserved under the action of $\pi(G)$, and in particular is normal in $K$.
 Since any element of $u^G$ has a finite orbit
under $\pi(G)$,  $N'$ contains only finite $\pi(G)$-orbits.
 In particular $N'$ is a finitely generated $FC$-group. It follows that $Tor(N')$ is a finite characteristic
subgroup of $N'$ and $N'/Tor(N')$ is free
 abelian with finite rank (cf. \cite{fc}).
Then one obtains a normal subgroup $N$ of $K$ satisfying condition $(i)$ by : if $Tor(N')\not=1$ then $N=Tor(N')$ and
otherwise $N=N'=\Z^n$.\hfill$\square$

\begin{fact} If $G\setminus K$ contains  $u^G$ finite, then either condition $(i)$ or $(ii)$ is
satisfied.\end{fact}

\noindent{\it Proof of the fact 4.}
 Let $u=k^{-1}q$ for some $k\in K$ and
$q\not=1$ lying in $Q$, such that $Z_G(u)$ has a finite index in $G$. Necessarily $q$ lies in $FC(Q)$, for $q^Q$ is the
image of $u^G$ under the projection of $G$ onto $Q$.

Let $h\in K$ and $\w=[u,h]\in K$ ; both $Z_G(u)$ and $Z_G(hu^{-1}h^{-1})$ have a finite index in $G$ and their
intersection lies in $Z_G(\w)$, so that $\w$ is an element of $K$ having a finite conjugacy class in $G$. If
$\w\not=1$, it follows from the fact 3  that condition $(i)$ is satisfied. So we suppose in the following that for any
$h\in K$, $[u,h]=1$, so that  $\pi_u$ is the identity on $K$. Hence $\theta_q$ is inner, for any $x\in K$,
$\theta_q(x)=x^k$.

Now let $Q_0=Z_G(u)\cap Z_Q(q)$, $Q_0$ is obviously contained in $Z_Q(k)$, so that $\theta(Q_0)$ is contained in
$Stab_\theta(k)$. Since $Q_0$ has a finite index in $Z_Q(q)$, it also has a finite index in $Q$, and then
$Stab_\theta(k)$ has a finite index in $\theta(Q)$, so that
 $k$ has a finite $\theta(Q)$-orbit. Hence condition $(ii)$ is satisfied.\hfill $\square$\\

\section{Formulation by mean of finite $\theta(Q)$-orbits}

One can formulate the theorem \ref{semidirect} by mean of the finite $\theta(Q)$-orbits in $K$.
\begin{thm}\label{thm2}
Let $G=K\rtimes_\theta Q\not=1$ and $O_\theta$ be the union of all finite $\theta(Q)$-orbits in $K$. Then $G$ is icc if
and only if :\begin{itemize}\item[$(a)$] $O_\theta\cap FC(K)=1$, and\item[$(b)$] $O_\theta\cap
\pi_{K}^{-1}(\theta(FC(Q)))=1$, and \item[$(c)$] the restricted homomorphism $\theta : FC(Q)\longrightarrow Aut(K)$ is
injective.
\end{itemize}
\end{thm}

\noindent{\it Proof.} Condition $(i)$ obviously implies that $O_\theta\cap FC(K)\not=1$. The converse is also true.
For, since $FC(K)$ is a characteristic subgroup of $K$, $O_\theta\cap FC(K)\not=1$ implies that $FC(K)$ contains a non
trivial finite $\theta(Q)$-orbit $O\not=\{1\}$. The union of conjugates of $O$ in $K$ is finite and preserved under
$\pi(G)$. So that for $k_0\in O$, $k_0^G$ is finite, $k_0\not=1$, and   condition $(i)$ follows from the fact 3 in the
proof of theorem \ref{semidirect}.

Conjunction of $(b)$ and $(c)$ is an immediate rephrasing
of the negation of condition $(ii)$. Conclusion follows from theorem \ref{semidirect}.\hfill$\square$\medskip\\
In particular, when $O_\theta=1$ one obtains a very concise statement.

\begin{cor}\label{cor1}
Let $G=K\rtimes_\theta Q\not=1$ such that all $\theta(Q)$-orbits in $K\setminus 1$ are infinite. Then $G$ is icc if and
only if the restricted homomorphism $\theta : FC(Q)\longrightarrow Aut(K)$ is injective.
\end{cor}

\section{On weakening condition $(ii)$} As we just have seen, in specific cases, condition $(ii)$ in theorem \ref{semidirect} can
be changed into the obvious : {\it $\theta : FC(Q)\longrightarrow Aut(K)$ is non injective}. Further examples follow
from :\pagebreak
\begin{prop}\label{prop1}
In the assumption of theorem \ref{semidirect}, if one moreover suppose at least one of the following conditions :
\begin{itemize}
 \item[--]   $K$ is abelian, \item[--] $K\setminus 1$ contains only infinite $\theta(Q)$-orbits,
 \item[--] the $\theta(Q)$-extension $\pi(G)$ of $Inn(K)$ splits, {i.e.} $\pi(G)=Inn(K)\rtimes \theta(Q)$,
\end{itemize} then condition (ii) can be strenghtened into :
\begin{itemize}
\item[--] the restricted homomorphism $\theta :FC(Q)\longrightarrow Aut(K)$ is non injective.
\end{itemize}
\end{prop}

\noindent{Proof.} If either $K$ is abelian or $\pi(G)=Inn(K)\rtimes \theta(Q)$, then necessarily one has that
$\theta(Q)\cap Inn(K)=1$ so that $\pi^{-1}_{K}(\theta(FC(Q)))=1$, and condition $(ii)$ becomes equivalent with $\theta
:FC(Q)\longrightarrow Aut(K)$ is non injective.\hfill$\square$\smallskip

In general one cannot strenghten condition $(ii)$ so far. For example if $K$ is icc and   $G=K\rtimes_\theta \Z$ then
$G$ is not icc each time  $\Theta:\Z\longrightarrow Out(K)$ is non injective ; which may happen while $\theta$ is
injective.

One may expect to weaken condition $(ii)$ into the condition that $\Theta:FC(Q)\longrightarrow Out(K)$ is non injective
; that is forgetting about hypothesis that $k$ has a finite $\theta(Q)$-orbit.
\begin{prop} \label{C1}
 In the assumption of theorem \ref{semidirect}, if one moreover suppose at least one of the
following conditions :
\begin{itemize}
 \item[--]   $Z(K)=1$, \item[--] $Q$ is finite or cyclic,
 \end{itemize} then condition $(ii)$ can be weakened into condition $(ii')$ :
\begin{itemize}
\item[$(ii')$]  $\Theta:FC(Q)\longrightarrow Out(K)$   is non injective.
\end{itemize}
\end{prop}

\noindent{\it Proof.} Obviously condition $(ii)$ implies condition $(ii')$. We prove the converse.\smallskip\\
{\it -- $Z(K)=1$.} Condition $(ii')$ implies that there exists $q\not=1$ in $FC(Q)$ such that $\theta_q$ is inner,
$\theta_q(x)=x^k$.
Any element $p\in Z_Q(q)$ is such that $\theta_p(k)\in k.Z(K)$. For, $\forall\, x\in K$,
$\theta_q(x)=k^{-1}xk=\theta_p\circ\theta_q\circ\theta_p^{-1}(x)=\theta_p(k^{-1})\,x\,\theta_p(k)$, implies that
$\theta_p(k)k^{-1}\in Z(K)$. So that with $Z(K)=1$, necessarily $\theta(Z_Q(q))$ lies in $Stab_\theta(k)$. Since
$Z_Q(q)$ has a
finite index in $Q$, $Stab_\theta(k)$ has a finite index in $\theta(Q)$, so that condition $(ii)$ is satisfied.\smallskip\\
 {\it -- $Q$ is finite or cyclic.} Suppose there exists $q\not=1$ in $FC(Q)$ such that $\theta_q$ is inner,
$\theta_q(x)=x^k$. Since $<q>_Q$ has a finite index in $Q$ and fixes $k$, condition $(ii)$ follows from condition
$(ii')$. (Moreover, if $Q$ is finite, condition $(i)$ is equivalent with $K$ not icc.)\hfill$\square$\medskip

We will see later several other particular cases for which the statement of theorem \ref{semidirect} becomes more
concise. But in general, condition $(ii)$ cannot be weakened into $(ii')$ as noted in the following remark.\smallskip\\
\noindent{\bf Remark 2.} Condition $(ii)$ of theorem \ref{semidirect} cannot in general be weakened in condition that
$\Theta:FC(Q)\longrightarrow Out(K)$ is non injective. For consider :
$$K=<a_1,a_2,k_1,k_2|\,[a_1,a_2],[a_i,k_j],i,j=1,2>\approx (\ZZ)\times F(2)$$
$$A=<a_1,a_2>_K\approx\ZZ\subset K$$
$$Q=<q_1,q_2|\,[q_1,q_2]>\approx\ZZ$$
Let $\theta_1\in Inn(K)$, {\it s.t.} $\forall\,x\in K,\theta_1(x)=x^{k_1}$ ; $\theta_1$ fixes $A$ pointwise. Let
$\theta_2\in Aut(K)$, {\it s.t.} $\theta_2$ is anosov on $A$, and $\theta_2(k_2)=k_2$, $\theta_2(k_1)=k_1\a$ for some
$\a\not=1$ lying in $A$. So defined, $\theta_1$ and $\theta_2$ commute, so that the map sending $q_1$ to $\theta_1$ and
$q_2$ to $\theta_2$ extends to an homomorphism $\theta : Q\longrightarrow Aut(K)$ ; moreover $\theta$ is injective.

Consider $G=K\rtimes_\theta Q$ ; we show that $G$ is icc despite that $\Theta:FC(Q)\longrightarrow Out(K)$ is non
injective. For any non trivial $x\in K$, $x^G$ is infinite, so that in particular condition $(i)$ of theorem
\ref{semidirect} is not satisfied. If condition $(ii)$ would be satisfied, it would follow that for some $n\geq 1$,
$k_1^n$ would have a $\theta(Q)$-finite orbit. We show that this cannot arise.

Consider $\theta_2\in\theta(Q)$, $\theta_2(k_1)=k_1\a$, $\a\not=1\in A$, so that for any $p\geq 1$,
$$\theta_2^p(k_1^n)=k_1^n\a^n\theta_2(\a^n)\theta_2^2(\a^n)\cdots\theta_2^{p-1}(\a^n)$$
Let $\phi_p:\ZZ\longrightarrow \ZZ$ be the map defined by $\phi_p(x)=x\theta_2(x)\theta_2^2(x)\cdots \theta_2^{p-1}(x)$
;
 $\phi_p$ turns out to be  an homomorphism. Let $M_\theta\in SL(2,\Z)$ be the matrix associated with $\theta_2$ ; it has
 two distinct irrational eigen values $\lambda_1,\,\lambda_2$. Let $M_{p}$ be the matrix associated with $\phi_p$.
 Then $M_{p}=Id+M_\theta+M^2_\theta+\cdots+M^{p-1}_\theta$. $M_{p}$ has two eigen values :
 $l_i=1+\lambda_i+\lambda_i^2+\cdots +\lambda_i^{p-1}$, $i=1,2$. They must be both non null because otherwise
 $\lambda_i^p=1$ which contradicts  that $M_\theta$ is anosov. Hence, for any $p\geq 1$,
  $\phi_p$ is injective. Since
for any $n\geq 1$,
 $\theta_2^p(k_1^n)=k_1^n\phi_p(\a^n)$, with $\a^n\not=1\in A$, the $\theta(Q)$-orbit of $k_1^n$
 is infinite, so that
 condition $(ii)$ is not satisfied. With theorem \ref{semidirect}, $G$ is icc, despite that the homomorphism
 $\Theta:FC(Q)\longrightarrow
 Out(K)$  is non injective.\\

\section{Further on weakening $(ii)$ : the twin $FC$-subfactors}
We keep on refining condition $(ii)$ by looking at what is in between the strenghtened condition
$\theta:FC(Q)\longrightarrow Aut(K)$ non injective, condition $(ii)$ and the weakened condition $(ii')$ :
$\Theta:FC(Q)\longrightarrow Out(K)$ for theorem \ref{semidirect}.\smallskip

\noindent{\bf Definition.} Let $G=K\rtimes_\theta Q$ ; We say that $C\rtimes C'$ is a {\it twin $FC$-subfactor of
$G$} when :\\
-- $C$ is a subgroup of $K$, $C\cap FC(K)=1$,\\
-- $C$ is $\theta(Q)$-stable with only finite $\theta(Q)$-orbits,\\
-- $C'$ is a normal subgroup of $Q$, $C'\subset FC(Q)$,\\
-- $\pi$ and $\theta$ are injective respectively on $C$ on $C'$ and $\pi(C)=\theta(C')$,\\
 (so that $\theta^{-1}\circ
\pi\,_{|C}:C\longrightarrow C'$ is an isomorphism),\\
-- $C\not=1$ (and so $C'$) is either finite or $\Z^n$.\smallskip\\
A twin $FC$ subfactor is either $\Z^n\times \Z^n$ or $C\rtimes_{Inn(C)} C$, for some  finite group $C$. It is
$\theta(Q)$-stable with only finite orbits, and $\theta^{-1}\circ \pi$ sends isomorphically the normal factor on the
retract one.
\begin{thm} Let $G=K\rtimes_\theta Q\not=1$ ; Then $G$ is not icc if and only if either condition (i) or at least one
of the following conditions is satisfied :
\begin{itemize}
\item[$(ii.a)$] $\theta:FC(Q)\longrightarrow Aut(K)$ is non injective, \item[$(ii.b)$] $G$ contains a twin
$FC$-subfactor.
\end{itemize}
\end{thm}

\noindent{\it Proof.} We first consider the sufficient part of the assumption. Condition $(i)$ implies that $G$ is not
icc follows from theorem \ref{semidirect} ; obviously condition $(ii.a)$ also implies $G$ not icc. If $G$ contains a
twin $FC$-subfactor $C\rtimes_{Inn(C)}C'$, then condition $(ii)$ of theorem \ref{semidirect} is satisfied with $q$
being any non trivial element of $C'$ and $k=\theta\circ\pi^{-1}(q)$, so that $G$ is not icc.\smallskip

We now prove the necessary part of the assumption. We suppose in the following that $G$ is not icc while it satisfies
neither condition $(i)$ nor condition $(ii.a)$ and prove that condition $(ii.b)$ must be satisfied.

With the theorem \ref{semidirect}, there exists $q\not=1$ in $FC(Q)$ and $k\not=1$ in K, such that $\theta_q(x)=x^k$
and $Stab_Q(k)$ has a finite index in $Q$. Let $C_Q$ be the subgroup of $FC(Q)$ finitely generated by $q^Q$ ; $C_Q$ is
a non trivial $FC$-group normal in $Q$.
Let $Q_1=\theta^{-1}(Stab_\theta(k))$, it has a finite index in $Q$. Clearly $Stab_\theta(k)$ is included in
$Z_{\theta(Q)}(\theta_q)$ ; if $Q_1\not\subset Z_Q(q)$, there would exist $p\in Q$ such that $[p,q]\not=1$ and
$\theta([p,q])=1$, which would contradict that $\theta:FC(Q)\longrightarrow Aut(K)$ is injective. Hence, $Q_1\subset
Z_Q(q)$, so that for $q_0=1,q_1,\ldots ,q_p$  a set of representatives of $Q/Q_1$,
 $C_Q$ is generated by the finite family $q, q^{q_1},\ldots , q^{q_p}$.

 Let $k_i$ be such that $k=\theta_{q_i}(k_i)$, then
$k_0=k,k_1,\ldots ,k_p$ is the $\theta(Q)$-orbit of $k$ and moreover $\theta_{q_i}^{-1}\circ
\theta_q\circ\theta_{q_i}(x)=x^{k_i}$. Let $C_K$ be the subgroup of $K$ generated by $k_0,k_1,\ldots ,k_p$ ; $C_K$ is
preserved under $\theta(Q)$ and  contains only finite $\theta(Q)$-orbits. An element in $C_K\cap FC(K)$ has a finite
conjugacy class in $G$, so that $C_K\cap FC(K)=1$ because otherwise as in fact 3 in the proof of theorem
\ref{semidirect} condition $(i)$ would follow.

By construction, $\pi(C_K)=\theta(C_Q)$. Each element of $Ker\, \pi_{|C_K}$ has a finite conjugacy class in $G$ so that
$\pi$ must be injective  on $C_K$ because otherwise, as in fact 3 in the proof of theorem \ref{semidirect}, condition
$(i)$ would follow. Moreover $\theta$ is injective on $C_Q$ because otherwise $\theta :FC(Q)\longrightarrow Aut(K)$
would be non injective. Hence $\theta^{-1}\circ\pi\,_{C_K}:C_K\longrightarrow C_Q$ is an isomorphism. Now $C_Q\not=1$
is a f.g. $FC$-group and hence with \cite{fc} either $Tor(C_Q)\not=1$ is a finite normal subgroup in $Q$, in which case
let $C'=Tor(C_Q)$, or $C_Q\approx \Z^n$, in which case let $C'=C_Q$. If $C$ denotes $\pi^{-1}\circ\theta(C')$, then
$C\rtimes_\theta C'$ is a twin $FC$-subfactor in $G$. \hfill$\square$\medskip

In conclusion suppose that $G=K\rtimes_\theta Q$ does not satisfy condition $(i)$ of theorem \ref{semidirect}. If
$\theta :FC(Q)\longrightarrow Aut(K)$ is non injective then $G$ is not icc. If $G$ is not icc despite $\theta$ is
injective then $G$ contains a twin $FC$-subfactor. It follows that $\Theta:FC(Q)\longrightarrow Out(K)$ is non
injective. If $\Theta$ is non injective, $G$ may be icc as seen in remark 2 ; $G$ is not icc whenever $G$ contains a
twin $FC$-subfactor.\smallskip

\noindent{\bf Example.} As in remark 2, consider :
$$K=<a_1,a_2,k_1,k_2|\,[a_1,a_2],[a_i,k_j],i,j=1,2>\approx (\ZZ)\times F(2)$$
$$A=<a_1,a_2>_K\approx\ZZ\subset K$$
$$Q=<q_1,q_2|\,[q_1,q_2]>\approx\ZZ$$
Let $\theta_1\in Inn(K)$, {\it s.t.} $\forall\,x\in K,\,\theta_1(x)=x^{k_1}$. Let $\theta_2\in Aut(K)$, {\it s.t.}
$\theta_2$ is anosov on $A$, and $\theta_2$ is the identity on $<k_1,k_2>_K$. Hence $\theta_1$ and $\theta_2$ commute
so that the map sending $q_1$ to $\theta_1$ and $q_2$ to $\theta_2$ extends to an injective homomorphism
$\theta:Q\longrightarrow Aut(K)$. So defined, $\theta(Q)$ fixes $k_1$ so that condition $(ii)$ of theorem
\ref{semidirect} is satisfied and $G$ is not icc ; condition $(i)$ is not satisfied since $FC(K)=A$ contains only
infinite $\theta(Q)$-orbits.
Let $C\subset K$, $C'\subset Q$ be generated respectively by  $k_1$ and $q_1$, then $C\rtimes_\theta C'=\Z\times\Z$ is
a twin $FC$-subfactor in $G$.

\begin{cor}
Let $G=K\rtimes_\theta Q$, {\it s.t.} $K$ does not contain any $\theta(Q)$-invariant subgroup $H\not=1$, either finite
or $\Z^n$ with only finite $\theta(Q)$-orbits.
Then $G$ is icc if and only if the restricted homomorphism $\theta :FC(Q)\longrightarrow Aut(K)$ is injective.
\end{cor}

\noindent{\it Proof.} Under these hypothesis, one cannot verifies condition $(i)$, and $G$ cannot contain any twin
$FC$-subfactor. \hfill$\square$

\section{Split extension of icc groups} We now consider the special case where at least one of the factors
is icc.

\begin{thm}\label{Kcci}
Let $G=K\rtimes_\theta Q$, with $K$ icc. Then $G$ is icc if and only if $\Theta:FC(Q)\longrightarrow Out(K)$ is
injective.
\end{thm}

\noindent{\it Proof.} Since $K$ is icc, on the one hand condition $(i)$ cannot arise and on the other $Z(K)=1$ so that
the conclusion follows from proposition \ref{C1} and theorem \ref{semidirect}.\hfill$\square$\smallskip

The following follows directly from theorem \ref{thm2}.
\begin{thm}\label{Qcci}
Let $G=K\rtimes_\theta Q$ with $Q$ icc. Then $G$ is icc if and only if the action of $\theta(Q)$ does not have any
finite orbit in $FC(K)\setminus 1$.
\end{thm}

\begin{cor}
The icc property is stable under split extension.
\end{cor}

\section{A particular case : word hyperbolic normal factor}

 Note  that a non elementary word
hyperbolic group $K$ is icc if and only if it does not contain any non trivial finite normal subgroup. For,
 in a word hyperbolic group the centralizer of any $\Z$ subgroup is virtually $\Z$. So that $FC(K)$ is a torsion
 $FC$-group. Since in a word hyperbolic group there is only finitely many conjugacy classes of torsion elements, $FC(K)$ is a
 finite normal subgroup of $K$.

\begin{thm} Let $G=K\rtimes_\theta Q$ with $K$ a non trivial word hyperbolic group.
Then $G$ is icc if and only if $K$ is icc and $\Theta:FC(Q)\longrightarrow Out(K)$ is injective.
\end{thm}

\noindent{\it Proof.} Note that, as already stated, $G=K\rtimes_\theta Q$ is not icc whenever $K\not=1$ is elementary.
So we suppose in the following that $G$ is non elementary. It follows that $Z(K)$ is finite
 (cf.
\cite{gro}). If $Z(K)\not=1$ then $K$ contains a non trivial finite normal subgroup so that both $K$ and $G$ are not
icc. If $Z(K)=1$, it follows from proposition \ref{C1} that $G$ is not icc if and only if it satisfies either condition
$(i)$ or condition $(ii')$. In a word hyperbolic group the normalizer of a $\Z$ subgroup is virtually $\Z$ ; in
particular, a non elementary hyperbolic group cannot contain any $\ZZ$ nor normal $\Z$ subgroup. Hence condition $(i)$
becomes that $K$ contains a non trivial finite normal subgroup, that is,  that $K$ is not icc.\hfill$\square$

\begin{cor} Let $G=K\rtimes_\theta Q$ with $K$ a non cyclic torsion free word hyperbolic group.
Then $G$ is icc if and only if $\Theta:FC(Q)\longrightarrow Out(K)$ is injective.
\end{cor}


\begin{thebibliography}{mot}

\bibitem[Co]{ydc}
Y.\textsc{de Cornulier}, \emph{Infinite conjugacy classes in groups acting on trees}, preprint (2005).


\bibitem[HP]{aogf3v}
P.\textsc{de la Harpe} et J.-P.\textsc{Préaux}, \emph{Groupes fondamentaux des variétés de dimension 3 et algèbres
d'opérateurs}, preprint arXiv:math.GR/0509449 v1 (2005).

\bibitem[Gr]{gro}
M.{\textsc{Gromov}}, \emph{Hyperbolic groups} in "\emph{Essays in group theory}", MSRI Publications, Springer (1987),
75--263.


\bibitem[ROIV]{roiv}
F.J.\textsc{Murray} et J.\textsc{von Neumann}, \emph{On rings of
operators}, IV, Annals of Math. {\bf 44} (1943), 716--808.

\bibitem[Ne]{fc}
W.\textsc{Neumann}, \emph{Groups with finite conjugacy classes},
Proc. of the L.M.S. {\bf 1} (1951), 178--187.

\bibitem[P1]{p1}
J.-P.\textsc{Préaux}, \emph{Wreath product of groups with infinite conjugacy classes}, preprint arXiV:math.GR/0612685
(2006).

\bibitem[P2]{p2}
J.-P.\textsc{Préaux}, \emph{Finite extension of group with infinite conjugacy classes}, preprint arXiV:math.GR/0703314
(2007).

\end{thebibliography}
 \end{document}